\documentclass{article}
\usepackage[a4paper, total={6in, 8in}]{geometry}
\usepackage{graphicx}
\usepackage{amsfonts,amsmath,amssymb,amsthm}
\usepackage{xcolor}
\usepackage{booktabs}
\usepackage{hyperref}
\usepackage{minted}

\let\phi\varphi

\newcommand{\prtt}[1]{\left(#1\right)}

\NewDocumentCommand{\tens}{t_}
 {%
  \IfBooleanTF{#1}
   {\tensop}
   {\otimes}%
 }
\NewDocumentCommand{\tensop}{m}
 {%
  \mathbin{\mathop{\otimes}\displaylimits_{#1}}%
 }

\newtheorem{defi}{Definition}
\numberwithin{defi}{section}
\newtheorem{exemp}[defi]{Example}
\newtheorem{prop}[defi]{Proposition}
\newtheorem{teo}[defi]{Theorem}

\hypersetup{
    colorlinks = true,
    linkcolor  = blue!70!black,
    citecolor  = green!50!black,
    urlcolor   = blue
}

\title{\texttt{LipschitzSaturation}: A \textit{Macaulay2} Package for Computing
       Lipschitz Saturations of Modules and Toric Varieties}
\author{Guilherme Schultz Netto \and Thiago da Silva}
\date{July 2026}

\begin{document}

\maketitle

\begin{abstract}
We introduce \verb|LipschitzSaturation|, a package for the computer algebra system \textit{Macaulay2} that implements algorithms for computing Lipschitz saturations of modules and toric varieties. In the module setting, the package handles three distinct saturation notions for an $\mathcal{O}_{X}$-submodule $\mathcal{M}\subseteq\mathcal{O}_{X}^{p}$: the 1-, 2-, and 3-Lipschitz saturations $\mathcal{M}_{S_{1}}$, $\mathcal{M}_{S_{2}}$ and $\mathcal{M}_{S_{3}}$, together with the auxiliary construction of the double module $\mathcal{M}_{D}$. To bypass the computationally intractable multivariate calculations for $\mathcal{M}_{S_{1}}$, we implemented a curve-based membership test, achieving near-constant runtime on parametric families that cause the purely algebraic method to time out or exhaust memory. For Toric Singularities, we implemented a construction algorithm. The package is freely available and requires \textit{Macaulay2} version 1.22 or later.  
\end{abstract}

\section{Introduction}
\par The concept of Lipschitz saturation of complex analytic algebras was introduced by Pham and Teissier~\cite{pham}, who observed that the germs of Lipschitz meromorphic functions on a reduced complex analytic space form an intermediate ring lying strictly between the coordinate ring $A$ and its integral closure $\overline{A}$. In the hypersurface case, this intermediate ring coincides with the Zariski saturation~\cite{zariskigts}, thereby linking bi-Lipschitz equisingularity with the purely algebraic theory of normalization. The algebraic properties of this construction were subsequently systematized by Lipman~\cite{lipman1}, who introduced the notion of relative Lipschitz saturation and established its key categorical properties, including idempotency and the fact that it always produces a radicial algebra extension. A modern treatment of these results connected to the theory of radicial algebras can be found in \cite{Da_Silva2024-dk}.

In practice, verifying whether a given element belongs to the Lipschitz saturation requires performing complex ideal membership tests via Gröbner basis reduction. These steps are routine for small examples, but become prohibitively expensive as the degrees of the generators or the number of variables grow. The \verb|LipschitzSaturation| package encapsulates all of this overhead, exposing a clean, user-facing interface. 

This document is organized as follows. Section \ref{sec:mod} treats the module setting, covering the double module construction, the $\rho$ functor, and the three saturation notions implemented by \verb|isMS1Element|, \verb|isMS1ElementOnCurve|, \verb|isMS2Element|, and \verb|isMS3Element|. Section \ref{sec:toric} discusses the algorithm to Toric Singularities.  

\section{The Double Module and Lipschitz Saturation of Modules}
\label{sec:mod}
The following definitions and results are due to \cite{gaffney2024generic}; new developments on the theory, mainly on the projective analytic spectrum by the same authors can be seen in \cite{daSilva_Gaffney_2026} and \cite{gaffney2024integral}.
Let $X \subseteq \mathbb{C}^n$ be a reduced analytic space and let
$\mathcal{M} \subseteq \mathcal{O}_X^p$ be an $\mathcal{O}_X$-submodule
finitely generated by global sections. Denote by
$\pi_1, \pi_2 \colon X \times X \to X$ the two coordinate projections. 

\begin{defi}
Let $h \in \mathcal{O}_X^p$. The \emph{double} of $h$ is
\begin{align*}
    h_D := (h \circ \pi_1,\; h \circ \pi_2)
           \in \mathcal{O}^{2p}_{X \times X}.
\end{align*}
The \emph{double} of $\mathcal{M}$, denoted $\mathcal{M}_D$, is the
$\mathcal{O}_{X \times X}$-submodule of $\mathcal{O}^{2p}_{X \times X}$
generated by $\{h_D \mid h \in \mathcal{M}\}$.
\end{defi}

Computing $\mathcal{M}_D$ directly from all elements of $\mathcal{M}$ is
impractical; the following proposition reduces it to the generators.

\begin{prop}[\cite{gaffney2024generic}]\label{prop:doublegen}
Suppose $\mathcal{M}$ is generated by $\{h_1, \ldots, h_r\}$.
Then $\mathcal{M}_D$ is also generated by each of the following sets:
\begin{enumerate}
    \item $\mathcal{B} = \bigl\{(h_1)_D, \ldots, (h_r)_D\bigr\}
          \cup \bigl\{\prtt{0,\; (z_i \circ \pi_1 - z_i \circ \pi_2)\,
          (h_j \circ \pi_2)} \mid
          i \in \{1,\ldots,n\},\; j \in \{1,\ldots,r\}\bigr\}$.
    \item $\mathcal{B}' = \bigl\{(h_1)_D, \ldots, (h_r)_D\bigr\}
          \cup \bigl\{\prtt{(z_i \circ \pi_1 - z_i \circ \pi_2)\,
          (h_j \circ \pi_2),\; 0} \mid
          i \in \{1,\ldots,n\},\; j \in \{1,\ldots,r\}\bigr\}$.
    \item $\mathcal{B}'' = \bigl\{(h_1)_D, \ldots, (h_r)_D\bigr\}
          \cup \bigl\{(z_i h_j)_D \mid
          i \in \{1,\ldots,n\},\; j \in \{1,\ldots,r\}\bigr\}$.
\end{enumerate}
\end{prop}

The function \verb|getDoubleModule(gensM)| implements generator set
$\mathcal{B}$ from Proposition~\ref{prop:doublegen}(1). It takes the
$p \times r$ generator matrix of $\mathcal{M}$, constructs the tensor product
ring $R \otimes_k R$ with its canonical projections, and assembles the
$2p \times (r + nr)$ generator matrix of $\mathcal{M}_D$ iteratively.

\begin{exemp}
\begin{verbatim}
i2 : B = QQ[x,y]

o2 = B

o2 : PolynomialRing

i3 : M = matrix {{x,y}, {y, x^2}}

o3 = | x y  |
     | y x2 |

             2      2
o3 : Matrix B  <-- B

i4 : Md = getDoubleModule M

o4 = | x_0 y_0   0             0              0             0                 |
     | y_0 x_0^2 0             0              0             0                 |
     | x_1 y_1   x_0x_1-x_1^2  x_0y_1-x_1y_1  y_0x_1-x_1y_1 y_0y_1-y_1^2      |
     | y_1 x_1^2 x_0y_1-x_1y_1 x_0x_1^2-x_1^3 y_0y_1-y_1^2  y_0x_1^2-x_1^2y_1 |

                                4                         6
o4 : Matrix (QQ[x , y , x , y ])  <-- (QQ[x , y , x , y ])
                 0   0   1   1             0   0   1   1
\end{verbatim}
\end{exemp}

\subsection{The \texorpdfstring{$1$}{1}-Lipschitz Saturation}

\begin{defi}\label{def:S1}
The \emph{$1$-Lipschitz saturation} of $\mathcal{M}$ at $x \in X$ is
\begin{align*}
\bigl(\mathcal{M}_{S_1}\bigr)_x
  := \bigl\{h \in \mathcal{O}^p_{X,x} \mid
     h_D \in \overline{\mathcal{M}_D} \text{ at } (x, x)\bigr\}.
\end{align*}
\end{defi}

To convert the geometric membership problem into an algebraic one, the package
employs the $\rho$ functor from~\cite{gaffney2024generic}. The integral closure of a
module lives in a free ambient module $\mathcal{O}_X^p$, which makes direct
computation hostile. The $\rho$ functor translates modules into ideals in a
larger polynomial ring, where standard integral-closure algorithms apply.

\begin{defi}
Given a vector $h = (h_1, \ldots, h_p) \in \mathcal{O}_X^p$ and projective
coordinates $[T_1 : \cdots : T_p]$, define
\begin{align*}
    \rho(h) := h_1 T_1 + \cdots + h_p T_p,
    \qquad
    \rho(\mathcal{M}) := \bigl\langle \rho(h) \mid h \in \mathcal{M} \bigr\rangle.
\end{align*}
\end{defi}

The key theoretical justification is:

\begin{teo}[\cite{gaffney2024generic}]\label{teo_rho}
Let $h \in \mathcal{O}^p_{X,x}$. Then
\begin{align*}
    h \in \overline{\mathcal{M}}
    \iff
    \rho(h) \in \overline{\rho(\mathcal{M})}
\end{align*}
at every point $(x, [T_1 : \cdots : T_p]) \in V(\rho(\mathcal{M}))$.
\end{teo}

Applying Theorem~\ref{teo_rho} to the double module $\mathcal{M}_D \subseteq
\mathcal{O}^{2p}_{X \times X}$, the package introduces $2p$ auxiliary
projective variables $T_1, \ldots, T_{2p}$ and works in the coordinate ring
$S = (R \otimes_k R)[T_1, \ldots, T_{2p}]$.

\paragraph{Algebraic criterion via Northcott--Rees reductions.}
Computing the full integral closure of $\rho(\mathcal{M}_D)$ in a ring with
$2n + 2p$ variables frequently exhausts memory, since the standard algorithm
constructs the entire Rees algebra of the ideal. Because the goal is merely to
certify membership of a \emph{single} element, the package instead invokes the
Northcott--Rees reduction criterion~\cite{swanson}: an element $z$ belongs to
the integral closure of an ideal $I$ if and only if $I$ is a reduction of
$J = I + \langle z \rangle$, i.e., there exists an integer $k \geq 1$ such that
\begin{equation}\label{eq:NR}
    J^k = I \cdot J^{k-1}.
\end{equation}
This transforms the normalization problem into an ideal-containment check, which
\textit{Macaulay2} resolves efficiently via Gröbner basis arithmetic in $S$.

The function \verb|isMS1Element(gensM, h)| implements this strategy:
\begin{enumerate}
    \item Compute $\mathcal{M}_D$ and the double vector
          $h_D = (h \circ \pi_1,\; h \circ \pi_2)^T$.
    \item Introduce $2p$ projective variables $T_1, \ldots, T_{2p}$ and
          form $S = (R \otimes_k R)[T_1, \ldots, T_{2p}]$.
    \item Compute $I = \rho(\mathcal{M}_D)$ and $z = \rho(h_D)$ in $S$.
    \item Set $J = I + \langle z \rangle$ and test
          $J^k \subseteq I \cdot J^{k-1}$ for $k = 1, \ldots, \texttt{MaxPower}$
          (default $5$). A positive answer at step $k$ certifies
          $h \in \mathcal{M}_{S_1}$.
\end{enumerate}

\begin{exemp}\label{ex:2.22}
We verify Example~2.22 of~\cite{gaffney2024generic}:
\begin{verbatim}
i2 : B = QQ[x,y]

i3 : M = matrix{{x,0,y},{y,x,0}}

i4 : h = matrix {{x},{3*y}}

i5 : isMS1Element(M,h)
o5 = false
\end{verbatim}
Thus $h = (x, 3y)^T \notin \mathcal{M}_{S_1}$.
\end{exemp}

\paragraph{Geometric criterion via analytic curves.}
For families parametrized by a degree parameter, the algebraic criterion suffers
from exponential complexity: the degrees of the generators of $\mathcal{M}_D$
grow with the parameter, driving Buchberger's algorithm into a combinatorial
explosion in S-polynomial reductions. The function
\verb|isMS1ElementOnCurve(gensM, h, curveVals)| circumvents this by pulling the
entire computation back along a user-supplied analytic curve
$\varphi \colon (\mathbb{C}, 0) \to (X \times X, (x, x))$. The user provides
the list of $2n$ coordinate pullbacks in $\mathbb{Q}[t]$. The package maps
$\mathcal{M}_D$ and $h_D$ into the univariate ring $S = \mathbb{Q}[t,
T_1, \ldots, T_{2p}]$, filters out zero-divisors, and resolves the membership
test $\rho(h_D) \in \overline{\rho(\mathcal{M}_D)}$ via polynomial division.
In the univariate domain, this operation runs in near-constant time regardless
of the degree parameter.

A \verb|false| return certifies $h \notin \mathcal{M}_{S_1}$; a \verb|true|
return is consistent with membership along that particular curve.

\begin{exemp}
Consider the Fernandes--Ruas example (\cite{gaffney2024generic}). For the parametric family
$F_n(x, y, z) = \frac{1}{3}x^3-z^2xy^{3n-2}+y^{3n}$, define the module
$\mathcal{M} = \langle \partial_x F_n, \partial_y F_n \rangle$ and the
candidate $h = \partial_z F_n$. Considering $n=1$ and the curve $\varphi(\tau) =
(\tau^2, \tau^2, \tau, -\tau^2, \tau^2, \tau)$ one can obtain:
\begin{verbatim}
i2 : R = QQ[x, y, z]
i3 : M = matrix {{x^2 - z^2*y, -z^2*x + 3*y^2}}
i4 : h = matrix {{-2*z*x*y}}
i5 : Kt = QQ[t]; tVar = Kt_0
i6 : curveVals = {tVar^2, tVar^2, tVar, -tVar^2, tVar^2, tVar}
i7 : isMS1ElementOnCurve(M, h, curveVals)
o7 = false
\end{verbatim}
\end{exemp}

Table~\ref{tab:runtimes} benchmarks the two strategies on the Fernandes--Ruas
family (AMD Ryzen~7 5700X). The degree parameter $n$ directly governs the
polynomial degree of the generators; as $n$ grows, the algebraic method becomes
intractable while the curve-based method remains essentially constant-time.

\begin{table}[h]
    \centering
    \begin{tabular}{@{}ccc@{}}
    \toprule
    $n$ & \texttt{isMS1Element} & \texttt{isMS1ElementOnCurve} \\
    \midrule
    $1$   & $0.272\,\mathrm{s}$   & $0.167\,\mathrm{s}$ \\
    $2$   & $1.229\,\mathrm{s}$   & $0.178\,\mathrm{s}$ \\
    $3$   & $33.76\,\mathrm{s}$   & $0.112\,\mathrm{s}$ \\
    $4$   & $26.00\,\mathrm{s}$   & $0.184\,\mathrm{s}$ \\
    $5$   & $443.7\,\mathrm{s}$   & $0.113\,\mathrm{s}$ \\
    $8$   & $1118\,\mathrm{s}$    & $0.111\,\mathrm{s}$ \\
    $10$  & \textsc{timeout}      & $0.114\,\mathrm{s}$ \\
    $20$  & \textsc{timeout}      & $0.193\,\mathrm{s}$ \\
    $50$  & \textsc{timeout}      & $0.213\,\mathrm{s}$ \\
    $100$ & \textsc{timeout}      & $0.184\,\mathrm{s}$ \\
    \bottomrule
    \end{tabular}
    \caption{Runtime comparison of \texttt{isMS1Element} and
             \texttt{isMS1ElementOnCurve} on the Fernandes--Ruas family.
             Timeout is set at $1200\,\mathrm{s}$.}
    \label{tab:runtimes}
\end{table}
\newpage
\subsection{The \texorpdfstring{$2$}{2}-Lipschitz Saturation}

\begin{defi}\label{def:S2}
The \emph{$2$-Lipschitz saturation} of $\mathcal{M}$ at $x \in X$ is
\begin{align*}
\bigl(\mathcal{M}_{S_2}\bigr)_x
  := \bigl\{h \in \mathcal{O}^p_{X,x} \mid
     \psi \cdot h \in \overline{(\psi \cdot \mathcal{M})_S}
     \text{ at } x,\;
     \forall\, \psi \colon X \to \mathrm{Hom}(\mathbb{C}^p, \mathbb{C})
     \bigr\}.
\end{align*}
\end{defi}

The condition requires that contracting $h$ with every $\mathcal{O}_X$-linear
form $\psi$ yields an element in the (1-Lipschitz) saturation of the
contracted module $\psi \cdot \mathcal{M}$. The function
\verb|isMS2Element(gensM, h)| tests this criterion via the ideal of $2 \times 2$
minors for modules of rank $2$:
$h \in \mathcal{M}_{S_2}$ if and only if
$J_2(h, \mathcal{M}) \subseteq \overline{J_2(\mathcal{M})}$,
where $J_2$ denotes the ideal of $2 \times 2$ minors. The proposition below
shows the two saturations are genuinely distinct.

\begin{prop}[\cite{gaffney2024generic}]
$\mathcal{M}_{S_1} \subseteq \mathcal{M}_{S_2}$ at every $x \in X$.
\end{prop}

Recalling Example~\ref{ex:2.22}, we now check the $2$-Lipschitz saturation:
\begin{verbatim}
i5 : isMS1Element(M, h)
o5 = false

i6 : isMS2Element(M, h)
o6 = true
\end{verbatim}
Since $h \in \mathcal{M}_{S_2}$ but $h \notin \mathcal{M}_{S_1}$, the
inclusion $\mathcal{M}_{S_1} \subseteq \mathcal{M}_{S_2}$ is strict in
this example.

\subsection{The \texorpdfstring{$3$}{3}-Lipschitz Saturation}
\begin{defi}\label{def:S3}
The \emph{$3$-Lipschitz saturation} of $\mathcal{M}$ at $x \in X$ is
\begin{align*}
\bigl(\mathcal{M}_{S_3}\bigr)_x
  := \bigl\{h \in \mathcal{O}^p_{X,x};
     J_k(h, \mathcal{M}) \subseteq \prtt{J_k(\mathcal{M})}_S
     \text{ at } x \bigr\},
\end{align*}
where $J_pk\mathcal{M})$ denotes the ideal of maximal ($k \times k$) minors
of any matrix of generators of $\mathcal{M}$, and
$J_k(h, \mathcal{M}) := J_k([\mathcal{M} \mid h])$ is the corresponding
ideal after appending $h$ as an extra column.
\end{defi}
The definition works great to implement the respective function, \verb|isMS3Element(gensM, h)|. It works testing this criterion directly,
using \textit{Macaulay2}'s \verb|minors| and \verb|integralClosure| commands together
with \verb|isSubset|. One has 
\begin{prop}[\cite{gaffney2024generic}]
    Suppose that $\mathcal{M}$ has generic rank $k$ on each component of $X$. Then, $\mathcal{M}_{S_2} \subseteq \mathcal{M}_{S_3}$.
\end{prop}

In general, $\mathcal{M}_{S_3}$ is the ``biggest" of the three module
saturations implemented. Still recalling Example \ref{ex:2.22}, one can get:
\begin{verbatim}
i5 : isMS1Element(M, h)
o5 = false

i6 : isMS2Element(M, h)
o6 = true

i7 : isMS3Element(M, h)
o7 : true
\end{verbatim}

\section{Toric Singularities}\label{sec:toric}
Expanding the context to toric singularities, we aim to construct the generator of the saturated algebra $\Gamma_S$, given a generator matrix $\Gamma$ of the semigroup. This part of the package was based on \cite{bernard2026lipschitzsaturationtoricsingularities}; although this article provide a criteria for membership, we can also expand it to construct an algorithm.

To compute the Lipschitz saturation $\Gamma_s$, the algorithm must operate within a well-defined finite search space. The saturated elements lie natively within the normalization $\overline{\Gamma}$, defined as the set of integer points situated in the dual cone $\sigma^\vee$ that also belong to the lattice generated by $\Gamma$. To guarantee finite termination, the computational procedure constructs a finite lattice grid determined by the maximum coordinate values of the original generators of $\Gamma$. Furthermore, the cone facets, obtained via the normal fans of the original semigroup, rigidly define the positive search space. This initial construction allows the algorithm to immediately discard candidate points situated outside the relevant geometry.

Once the initial grid is generated, the algorithm proceeds to filter the normalization set. Every candidate $\gamma \in \overline{\Gamma}$ extracted from the bounding box is evaluated against three rigorous criteria:

\begin{enumerate}
    \item For every facet of the cone, the projection of the candidate $\gamma$ must not exceed the maximum projection yielded by the original generators, ensuring that $\gamma$ does not violate the upper bounds of the singularity's structure.
    \item To rigorously verify if $\gamma$ is positioned inside the Newton Polyhedron, the algorithm maps the problem to a homogenized cone construction in $\mathbb{Z}^{d+1}$ dimensions, eliminating elements that fall strictly below the characteristic boundary of the singularity.
    \item Bypassing standard scalar modulo operators, the package implements a rigorous module map $\phi: \mathbb{Z}^n \to \mathbb{Z}^d$. This homological constraint guarantees that the candidate $\gamma$ is a valid integer linear combination of the semigroup generators.
\end{enumerate}

The execution flow translates these theoretical constraints into a sequential reduction pipeline: it inputs the generator matrix, generates the normalization grid incorporating the bounding and cone constraints, applies the geometric filters concerning directional projections and the homogenized convex hull, evaluates the surviving elements against the image of the module map $\phi$, and outputs the saturated matrix. Evaluating elements against the image of $\phi$ bypasses the arithmetic failures inherent to standard divisibility checks, specifically in cases where the image does not span the full lattice $\mathbb{Z}^d$. Consequently, this multidimensional filtering approach scales efficiently into higher dimensions compared to naive exhaustive coordinate searches.

\begin{exemp}
To validate the computational implementation in higher dimensions, we consider a non-isolated 3D singularity constructed by extending the foundational 2D example from \cite{bernard2026lipschitzsaturationtoricsingularities}. We define the generator matrix in $\mathbb{Z}^3$, incorporating a gap that requires saturation:

\begin{verbatim}
i1 : gensMat = matrix {{6, 5, 0, 0}, {0, 3, 6, 0}, {0, 0, 0, 6}}

o1 = | 6 5 0 0 |
     | 0 3 6 0 |
     | 0 0 0 6 |

              3       4
o1 : Matrix ZZ  <-- ZZ
\end{verbatim}

The algorithm systematically sifts through the lattice generated by the bounding box. It rejects points that fall outside the homogenized convex hull or those that fail the strict lattice divisibility mapped by $\phi$. The procedure isolates the exact saturated monoid element $\left[ 3, 3, 0 \right]^T$, yielding the following matrix:

\begin{verbatim}
i2 : sat = getLipschitzSaturationToric(gensMat)

o2 = | 6 5 0 0 3 |
     | 0 3 6 0 3 |
     | 0 0 0 6 0 |

              3       5
o2 : Matrix ZZ  <-- ZZ
\end{verbatim}

\end{exemp}

The final generator matrix correctly appends the missing monomial, confirming that the implementation accurately recovers the theoretical expectations for arbitrary dimensions.
\bibliography{ref2}
\bibliographystyle{abbrv}

\end{document}